\documentclass[11pt]{amsart}

\usepackage{amsmath, amssymb, amsthm}
\usepackage[hidelinks]{hyperref}
\usepackage[parfill]{parskip}        % comment out to return to indent paragraphs.
\usepackage[margin=1.05in]{geometry} 
\usepackage{tikz}
\usetikzlibrary{decorations.markings}
\usepackage{tikz-cd}
\usetikzlibrary{matrix}
\usepackage{comment}
\usepackage{mathtools}

\title{Planar aperiodic tile sets: from Wang tiles to the Hat and Spectre monotiles}\keywords{} \subjclass{}

\author[Tinka Bruneau]{Tinka Bruneau}
 
\author[Michael F.\ Whittaker]{Michael F.\ Whittaker}
\address{School of Mathematics and Statistics, University of Glasgow, University Place, Glasgow Q12 8QQ, United Kingdom}
\email{2580916b@student.gla.ac.uk, Mike.Whittaker@glasgow.ac.uk}

\def\N{\mathbb{N}}

\def\PP{\mathcal{P}}

\def\ep{\varepsilon}

\newcommand{\interior}{\mathrm{interior}}

\newtheorem{theorem}{Theorem}[section]

\theoremstyle{definition}

\numberwithin{equation}{section} % Number equations by section
\numberwithin{figure}{section}   % Number figures by section

\tikzstyle{vertex}=[circle]
\tikzstyle{goto}=[->,shorten >=1pt,>=stealth,semithick]

\makeatletter
\@namedef{subjclassname@2020}{%
  \textup{2020} Mathematics Subject Classification}
\makeatother

\thanks{We'd like to thank Michael Baake, Kevin Brix, Felix Flicker, Robbert Fokkink, Franz G\"ahler, Craig Kaplan, Jan Maz\'{a}\v{c} and Jamie Walton for excellent comments and suggestions on early drafts. This paper was written while the second author was a guest at the Fields Institute, and he thanks them for their hospitality and exceptional research environment.}
\begin{document}
\maketitle

\begin{abstract} 
A brief history of planar aperiodic tile sets is presented, starting from the Domino Problem proposed by Hao Wang in 1961. We provide highlights that led to the discovery of the Taylor--Socolar aperiodic monotile in 2010 and the Hat and Spectre aperiodic monotiles in 2023. The Spectre tile is an amazingly simple monotile; a single tile whose translated and rotated copies tile the plane but only in a way that lacks any translational periodicity. We showcase this breakthrough discovery through the 60$+$ years that aperiodic tile sets have been considered.
\end{abstract}

\section{Introduction}

In 1961, Hao Wang asked if it is possible to algorithmically determine whether translated copies of a finite set of marked square tiles can tile the plane. This became known as the Domino Problem. As a first attack on the Domino Problem, Wang posed his \emph{Fundamental Conjecture}: a set of Wang tiles either tiles the plane with translational periodicity or does not tile. The other possibility is a set of tiles that can form a tiling of the plane but always without any translational periodicity, such a tile set is called \emph{aperiodic}.

Robert Berger found a counterexample to Wang's Fundamental Conjecture by finding an aperiodic tile set with 20,426 tiles. This proved undecidability of the Domino Problem. After finding such a set of tiles, the search was on for the minimal size of such a set. Emmanuel Jeandel and Micha\"el Rao recently proved that an aperiodic set of Wang tiles must have at least 11 tiles and gave an example of such a tile set. We provide pictures of all tilings mentioned in the introduction throughout the article and the reader is invited to look ahead to get a feeling for these tile sets and their tilings.

Wang tiles are always squares, but what if we allow tiles to come in arbitrary shapes and allow tiles to be rotated and reflected? In the 1970s, both Sir Roger Penrose and Robert Ammann found aperiodic tile sets with just two tiles, now commonly referred to as the Penrose and Ammann--Beenker tiles, respectively. This left open the aperiodic monotile problem: Is there a single tile that can tile the plane but only in a way that forbids translational periodicity?

Two early attempts found monotiles that use overlaps to force aperiodicity. While these are not our focus, they are very interesting both historically and from the standpoint of aperiodic order. The first of these was Gummelt's covering of the Penrose tiling by a single decorated decagon \cite{Gum}. Gummelt's overlapping tile was really the first mono-cluster, and provided the first evidence that a monotile could be possible. The second was Penrose's $(1+\ep+\ep^2)$-tiling \cite{BGG,Pen2}. This was originally described as a single hexagon with a smaller hexagon and half-hexagons inscribed in its interior. The matching rules required overlapping hexagons to form the same type of hexagons as at the larger scale. An excellent discussion of both of these tiling can be found in \cite{BG}.

In a breakthrough result of Joan Taylor in 2010 \cite{Tay}, the first connected and completely geometrically defined monotile was discovered, but the tile is not simply connected. That is, the tile is not a topological disc. Taylor joined forces with Joshua Socolar to introduce their tile to the mathematical community \cite{ST,ST2}. The original formulation of the Taylor--Socolar tile was described as a decorated hexagon with certain matching rules between both neighbouring and next to neighbouring tiles, and the tiling requires a reflected tile.

In March 2023, Dave Smith, Craig Kaplan, Joseph Myers, and Chaim Goodman-Strauss found a simply connected monotile that they called the Hat \cite{Hat}. Amazingly, the Hat is a simple polygonal shape and is entirely geometric, in the sense that no additional matching rules on how tiles are allowed to meet are required to enforce aperiodicity. What's more, in the article, a beautiful proof of aperiodicity was used to show that there are, in fact, an uncountable number of tiles in the Hat family that are also simply connected monotiles. The tiling requires a reflected version of the Hat. This left open the question of whether a simply connected and geometric solution to the monotile problem without reflections is possible.

In a marvellous stroke of insight, Smith, Kaplan, Myers, and Goodman-Strauss realised that a member of the Hat family could also be used to tile the plane without using a reflected copy of the tile \cite{Spectre}. The Spectre aperiodic monotile was discovered. The Spectre is in the family of Hat tilings, but was one of the singular members that allowed periodic tilings. The authors realised that forbidding the reflected tile still allowed tilings of the plane, but all such tilings lack translational periodicity. Thus, the Spectre provides a remarkable solution to the aperiodic monotile problem!

\section{Tiles, tilings, and their properties}

In this article, we restrict ourselves to two dimensional tilings. We mostly follow the language and terminology defined in Baake and Grimm's \emph{Aperiodic Order} \cite{BG}, which is highly recommended to the interested reader. The second main reference on general tiling theory is the book of Gr\"unbaum and Shephard \cite{GS}, which contains an excellent introduction to Wang tiles and aperiodic tile sets.

We begin with the building blocks of a tiling. A \emph{prototile} is a labelled subset of $\mathbb{R}^2$ that is equal to the closure of its interior. This just means that there are no dangling bits hanging from the tiles; in fact, one should probably just think of labelled polygons. In addition, we often allow \emph{decorations} of the prototiles, like edge markings that must match, arrows to determine tile orientation, or lines that must continue across tile edges. Often, but not always, these decorations can be realised by puzzle like bumps and dents in the tile edges. We'll see examples of this soon.

Let $\PP$ be a finite set of prototiles and $G$ a subgroup of the isometry group of $\mathbb{R}^2$. A \emph{tiling} of $\mathbb{R}^2$ (the plane) is a countable collection of \emph{tiles} $T=\{t_i \mid i\in\mathbb{N}\}$ such that:
\begin{itemize}
\item $t_i=\gamma \cdot p$ for some $p\in\PP$ and $\gamma$ in $G$;
\item $\displaystyle{\bigcup_{i\in\N}}\hspace{0.1cm}t_i=\mathbb{R}^2$;
\item $\interior{(t_i)}\cap\interior{(t_j)}=\varnothing$ if $i \neq j$.
\end{itemize}
Let's take a minute to understand the definition. The first bullet point specifies that all tiles in the tiling must be isometric copies of prototiles (i.e., some composition of a translation, rotation and/or reflection of a prototile), where the group $G$ specifies the exact subset of isometries we allow for a given prototile set. Here we are thinking of the prototiles as actually sitting in the plane, so that we can move them about to form a tiling. In this article, $G$ will always be the translation group, the direct isometry group (translations and rotations), or the full isometry group. Note that when we use the full isometry group $G$, we will include the reflected tile(s) in the prototile set for clarity.  The second bullet point implies that the tiling covers the entire Euclidean plane and the third specifies that tiles never overlap except at their edges. Moreover, in the case that we have extra decorations, the tiling must also satisfy any extra rules specified by the decorations.

\begin{figure}[ht]
\begin{center}
\includegraphics[width=15cm]{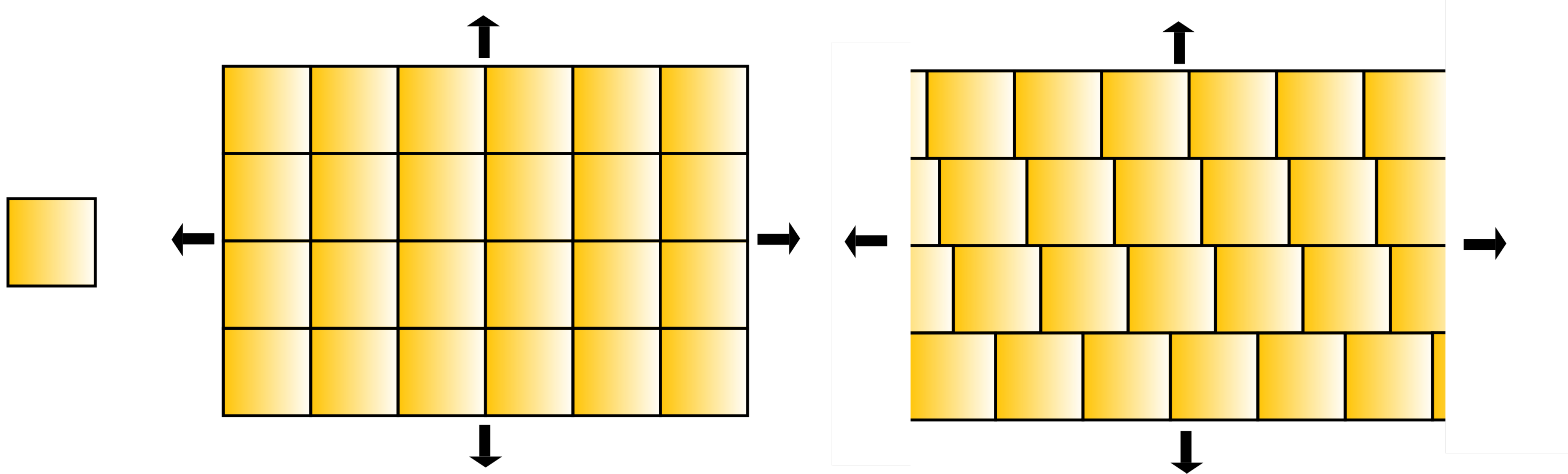}
\end{center}
\caption{Tilings of the plane by a single square prototile. The prototile is on the left along with two possible tilings.}
\label{SquareTiling}
\end{figure}

As Figure \ref{SquareTiling} depicts, it is easy to find tilings of the plane by a single square prototile of side length one. In fact, there are infinitely many possible tilings that arise from an unmarked square prototile. We note that it's impossible to pictorially represent a complete tiling of the plane, and the reader is meant to extrapolate how a complete tiling is produced from the small patch provided. For both examples in Figure \ref{SquareTiling}, we can let $G$ be the subgroup of translations. The tiling on the right reveals why there are infinitely many possible square tilings, even taken up to translation, by choosing different relative shifts for rows of square tiles. In this case, the second row is shifted by 1/2 with respect to the bottom row that's sitting on the $x$-axis, and the third is shifted by $1/3$ with respect to the bottom, and so on. We could use the same procedure to form rows below the bottom row to make a complete tiling of the plane.

We can translate a tiling $T$ by a vector $x$ in $\mathbb{R}^2$ via $T+x:=\{t + x \mid t \in T\}$. Notice that the tiling in the middle of Figure \ref{SquareTiling} has the property that $T+(1,0)=T$. Indeed, if we translate all tiles in the tiling over by one unit to the left the sets $T+(1,0)$ and $T$ are exactly the same! We say that a tiling $T$ is \emph{periodic} if there exists a nonzero translation $x$ in $\mathbb{R}^2$ such that $T+x=T$, and we say that $T$ is \emph{nonperiodic} if $T+x=T$ implies $x=0$. That is, if we took an infinite transparent photocopy of the tiling $T$, then the tiling is periodic if we can shift this photocopy by some non-trivial translation so that it perfectly matches $T$ and nonperiodic if it only matches in exactly one place. The tiling on the right of Figure \ref{SquareTiling} is also periodic with $x=(1,0)$. We'll see examples of nonperiodic tilings in the next section.

There are many more concepts that are important in tiling theory. Again, the book of Baake and Grimm \cite{BG} contains a modern treatment and an exceptional foreword written by Roger Penrose. For example, there has been discussion on the \emph{local indistinguishability classes} (LI-classes) for the Hat and Spectre monotiles, see \cite[p.2]{BGS}, where the notion of LI-classes can be found in \cite[Section 5.1.1]{BG}.

\section{Wang tiles}

A founder of modern tiling theory was the philosopher Hao Wang (1921--1995). A set of square prototiles with marked edges are called Wang tiles. Here we label the edges with a colour and also a number, but merely to aid legibility. Tiles must meet along complete edges only when the symbol matches and we only allow $G$ to be the group of translations; no rotations or reflections of prototiles allowed.

In 1961 Wang \cite{Wang} proposed the Domino Problem: given a set of Wang tiles, is it possible to algorithmically determine whether the set tiles the plane? Wang tiles are theoretically important in logic, since the behaviour of any Turing machine can be mimicked using some particular set of Wang tiles, see \cite[Section 11.4]{GS}.

\begin{figure}[h]
\begin{center}
\includegraphics[width=15cm]{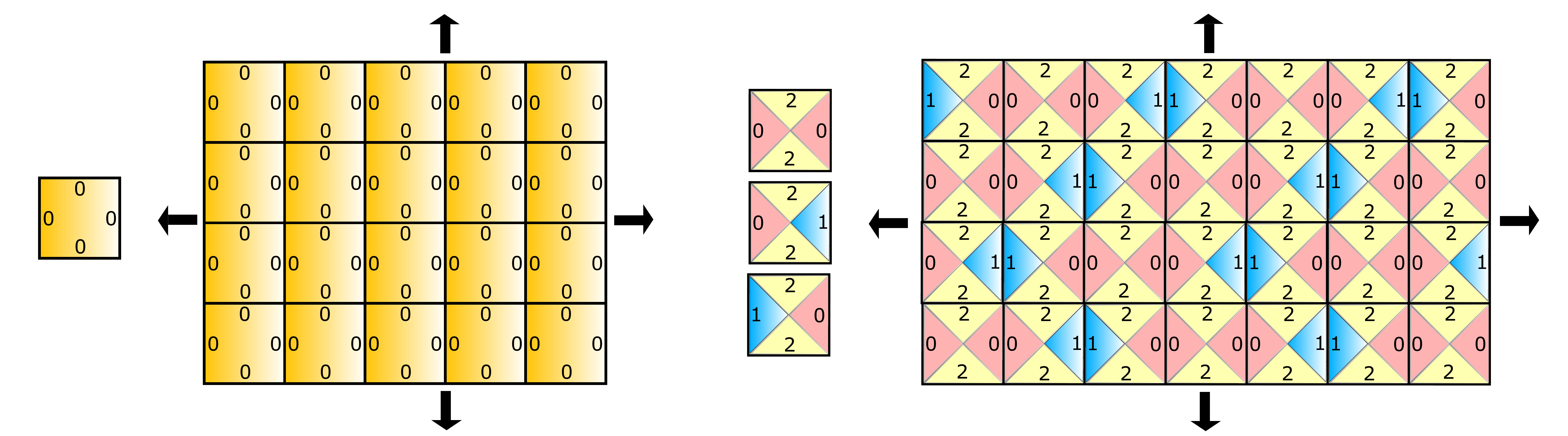} 
\end{center}
\caption{Two sets of Wang tiles, the left can only tile periodically while the right can tile both periodically and nonperiodically.}
\label{Wang_per}
\end{figure}

Wang realised that a set $\PP$ of Wang tiles has four possibilities:
\begin{enumerate}
\item The set $\PP$ does \textbf{not} tile the plane, such as a single tile with four distinct symbols on its edges;
\item The set $\PP$ can only tile the plane periodically, such as the Wang tile on the left of Figure \ref{Wang_per};
\item The set $\PP$ can tile both periodically and nonperiodically, such as the Wang tile set on the right of Figure \ref{Wang_per} (points for guessing a method to fill the tiling out to the plane in a non-periodic way);
\item The set $\PP$ tiles the plane but only nonperiodically.
\end{enumerate}
Wang realised that the existence of a Wang tile set satisfying (4) implies that the Domino Problem is undecidable \cite{Wang}. This led to his \emph{Fundamental Conjecture} that there are no Wang tile sets satisfying (4). However, Wang's student, Robert Berger, found the first such Wang tile set with 20,426 tiles \cite{Ber}! According to Gr\"unbaum and Shephard \cite[Chapter 11]{GS}, Berger subsequently got the set down to only 104 tiles (actually 103 tiles, see \cite[Section 1.2]{JR}). Shortly after, Donald Knuth (the inventor of TeX amongst other things) modified Berger's set to only require 92 tiles. The hunt was on to find the smallest set of Wang tiles, with contributions from Hans L\"auchli, Raphael Robinson, Roger Penrose, Robert Ammann, and John Conway finding Wang tiles sets for $n=56, 52, 40, 35, 34, 32, 24, \text{ and } 16$, see \cite[Section 11.1]{GS}. The following quote appears in \cite[p.596]{GS}:
\begin{quote}
The reduction in the number of Wang tiles in an aperiodic set from over 20,000 to 16 has been a notable achievement. Perhaps the minimum possible number has now been reached. If, however, further reductions are possible then it seems certain that new ideas and methods will be required.
\end{quote}

\begin{figure}[ht]
\begin{center}
\includegraphics[width=10cm]{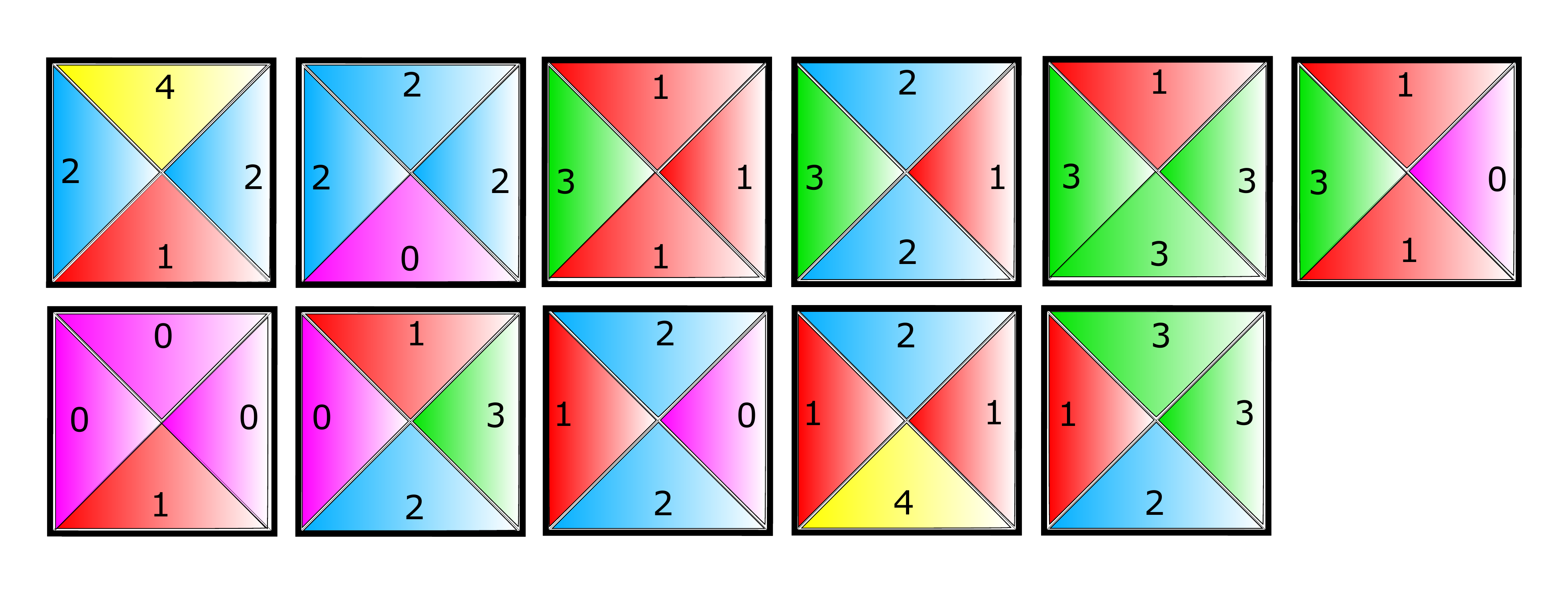} 
\end{center}
\caption{The smallest possible set of Wang tiles was found by Jeandel and Rao \cite{JR} consisting of just 11 tiles.}
\label{JR_tiles}
\end{figure}

\begin{figure}[ht]
\begin{center}
\includegraphics[width=15cm]{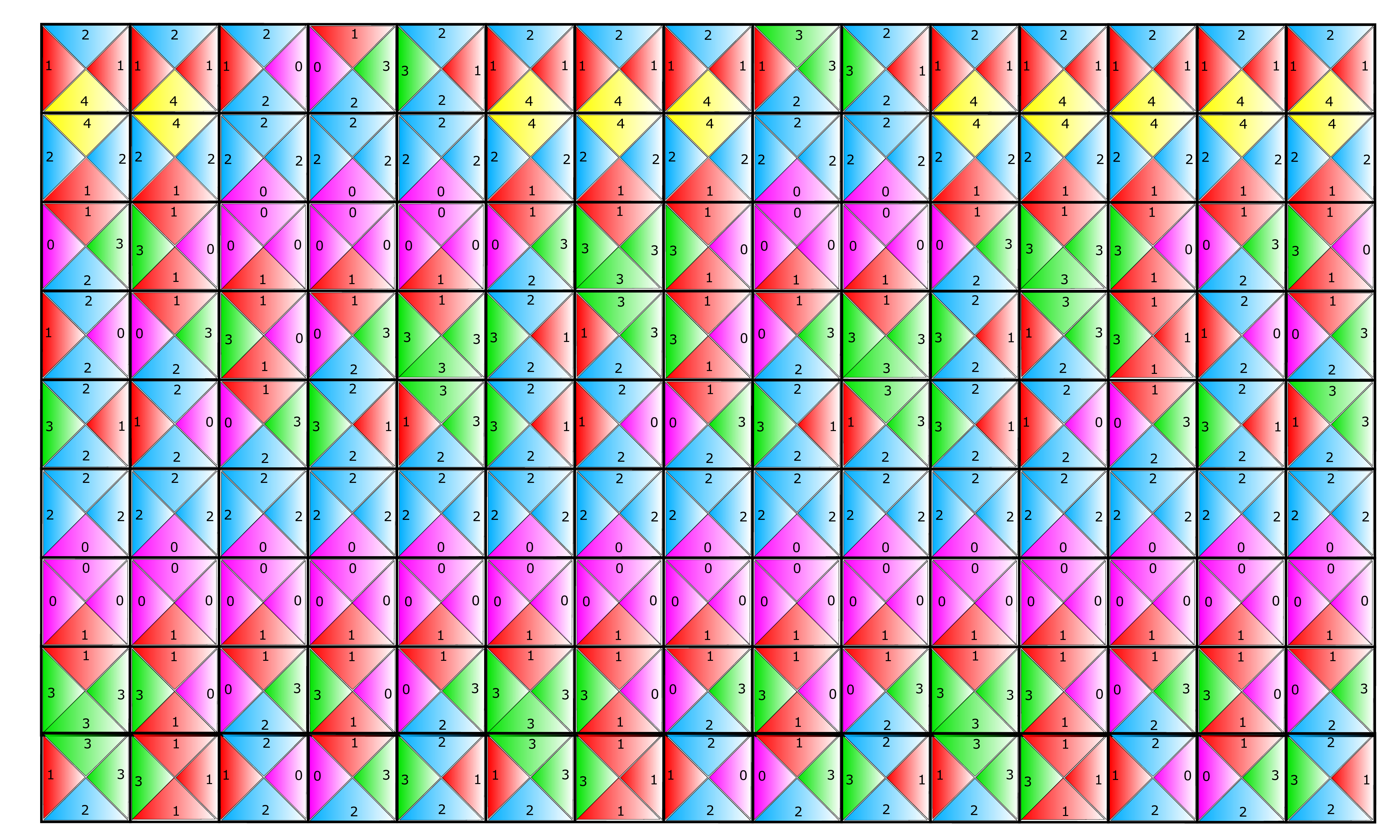} 
\end{center}
\caption{A patch of Jeandel--Rao Wang tiles shows the type of complexity required to reduce the number of tiles to 11.}
\label{JR_patch}
\end{figure}

Two further reductions were found, the Kari--Culik Wang tile set with only 13 tiles was found using number theoretic methods, see \cite[Section 5.7.4]{BG}. The Jeandel--Rao Wang tile set in Figure \ref{JR_tiles} reduced this number to just 11 Wang tiles \cite{JR}! The new method in this case was an exhaustive computer search, so we now know that this is the minimal number, although we don't know that the depicted Wang tile set is essentially unique as an aperiodic set of 11 Wang tiles. A patch of Jeandel--Rao tiles is in Figure \ref{JR_patch}.

\section{Relaxing the rules and the Penrose tiles}

A set of prototiles, of any shape, that admits tilings of the plane but only nonperiodically is called an \emph{aperiodic tile set}. During the flurry of activity to reduce the number of Wang tiles in the 1970s, Roger Penrose and Robert Ammann began considering aperiodic tile sets of polygons with specified markings. Both independently found an aperiodic tile set with just two tiles! However, we focus only on the Penrose tiles here. An excellent introduction to aperiodic tiles sets, with proofs that had not appeared previously in the literature, is \cite[Chapter 10]{GS}.

\begin{figure}[ht]
\begin{center}
\includegraphics[width=10cm]{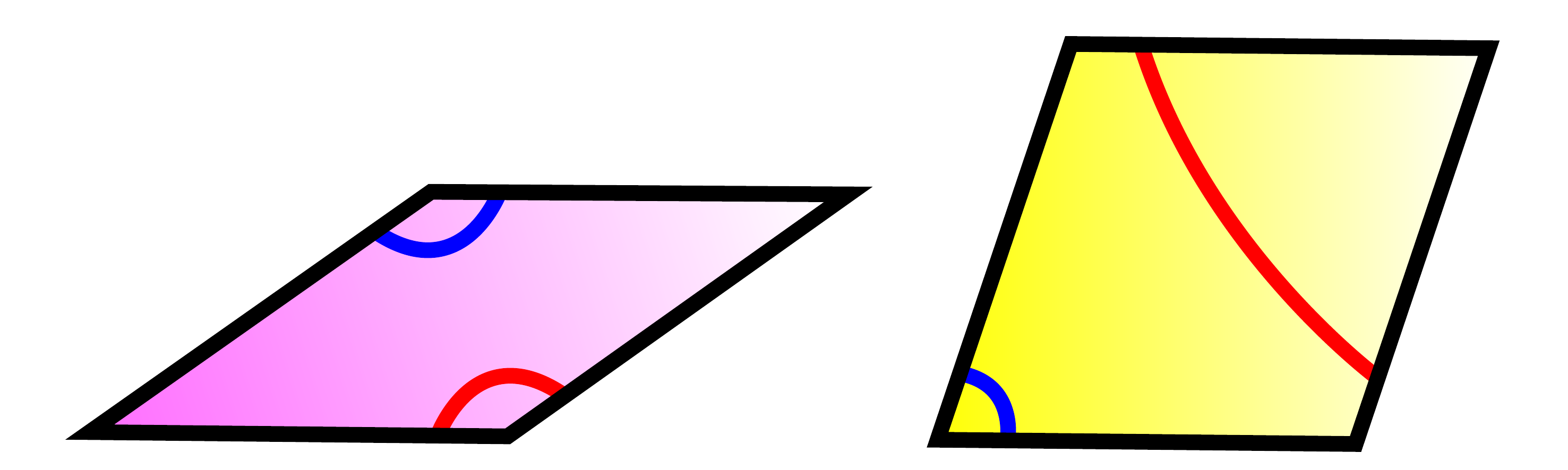} 
\end{center}
\caption{The two tile Penrose aperiodic tile set. The tiles must meet full-edge to full-edge and the red and blue lines must continue from one tile to the next.}
\label{Penrose_tiles}
\end{figure}

\begin{figure}[ht]
\begin{center}
\includegraphics[width=15cm]{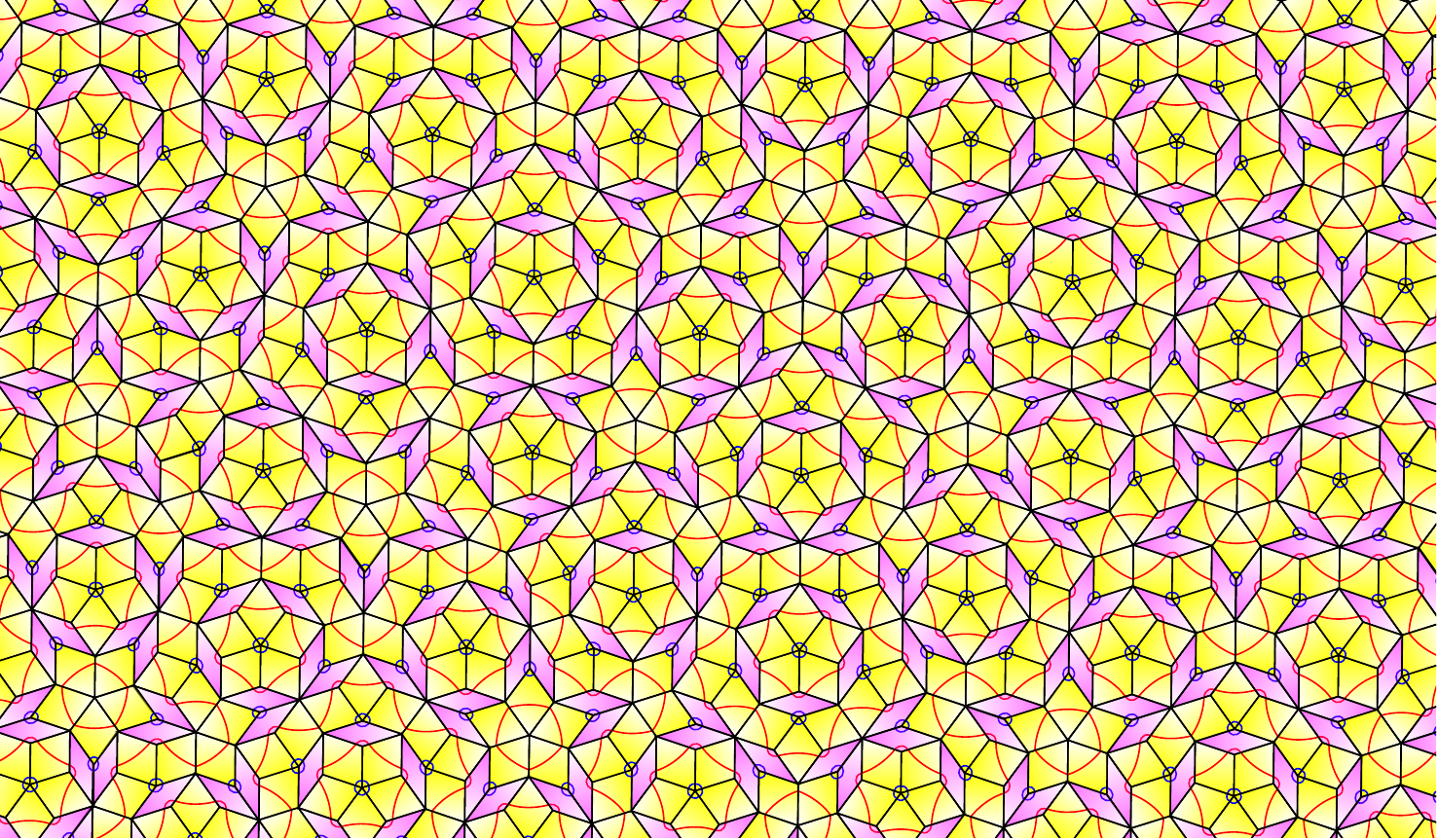} 
\end{center}
\caption{A patch of a Penrose tiling.}
\label{Penrose_patch}
\end{figure}

The Penrose two tile aperiodic set is depicted in Figure \ref{Penrose_tiles} and a patch of tiles in Figure \ref{Penrose_patch}. The tiles are simple polygonal shapes that must match full-edge to full-edge with the blue and red lines continuing from one tile to the next. We note that the line matching rules can be encoded purely geometrically by puzzle like bumps and dents, see \cite[p.154]{BG}. These tiles were reduced from other versions of Penrose tile sets, with more prototiles, and have a 10 fold symmetry group.

\begin{figure}[ht]
\begin{center}
\includegraphics[width=7cm,angle=-90]{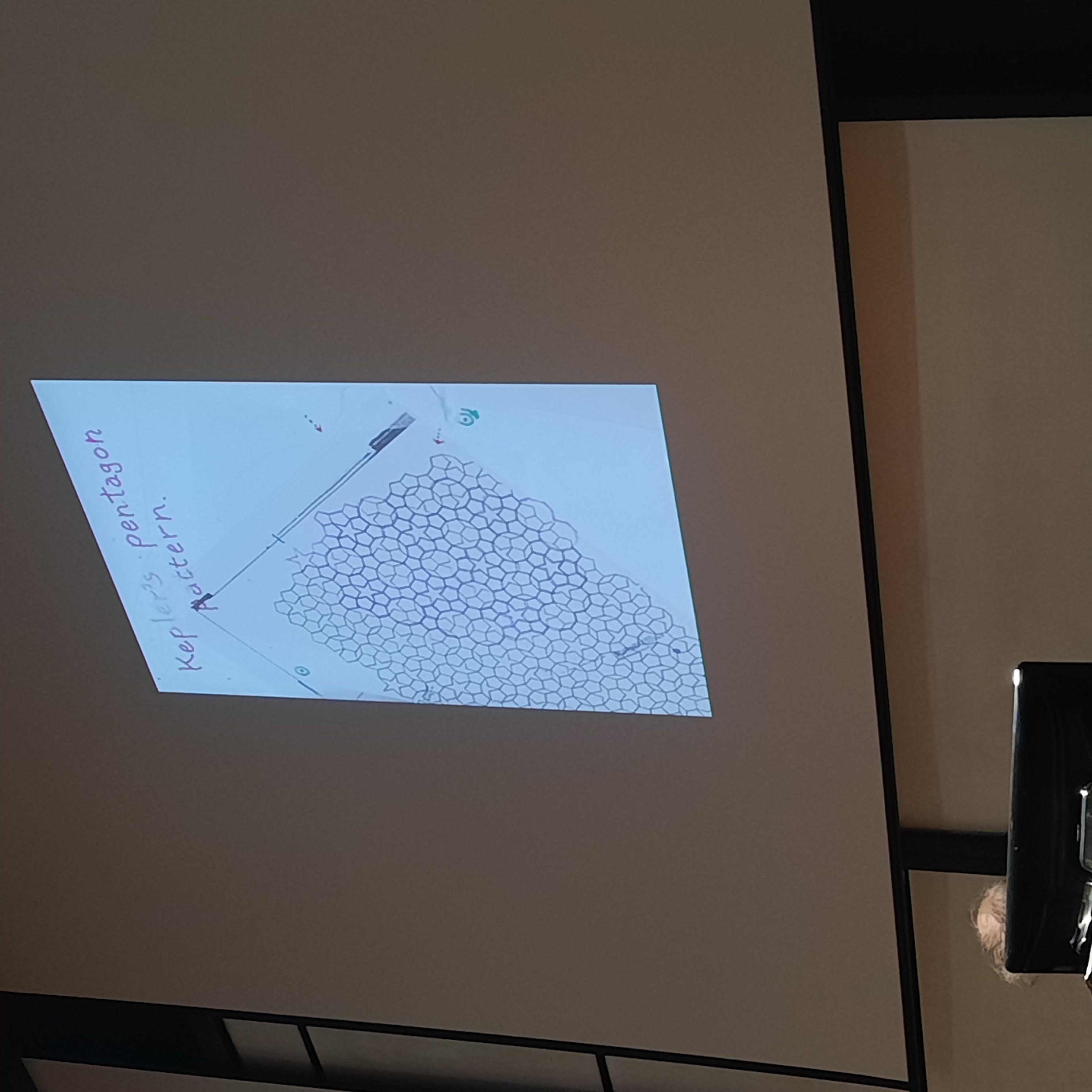} 
\end{center}
\caption{Roger Penrose overlaying a patch of the Penrose tiling over Kepler's pentagon pattern.}
\label{Penrose_kepler}
\end{figure}

Roger Penrose gave a lecture at \emph{Hatfest: celebrating the discovery of an aperiodic monotile} at the University of Oxford in July, 2023 where he showed how he constructed his two tile set in Figure \ref{Penrose_tiles}. He suggested that he may have been indirectly inspired by a book of Kepler that was on his parents' bookshelf. In fact, during his lecture he overlayed his tiles with an image produced by Kepler to an almost perfect match, see Figure \ref{Penrose_kepler}. More information on this link can be found in Rodrigo Trevi\~no's Notices of the American Mathematical Society article \cite{Trev}.

\begin{figure}[ht]
\begin{center}
\includegraphics[width=8cm]{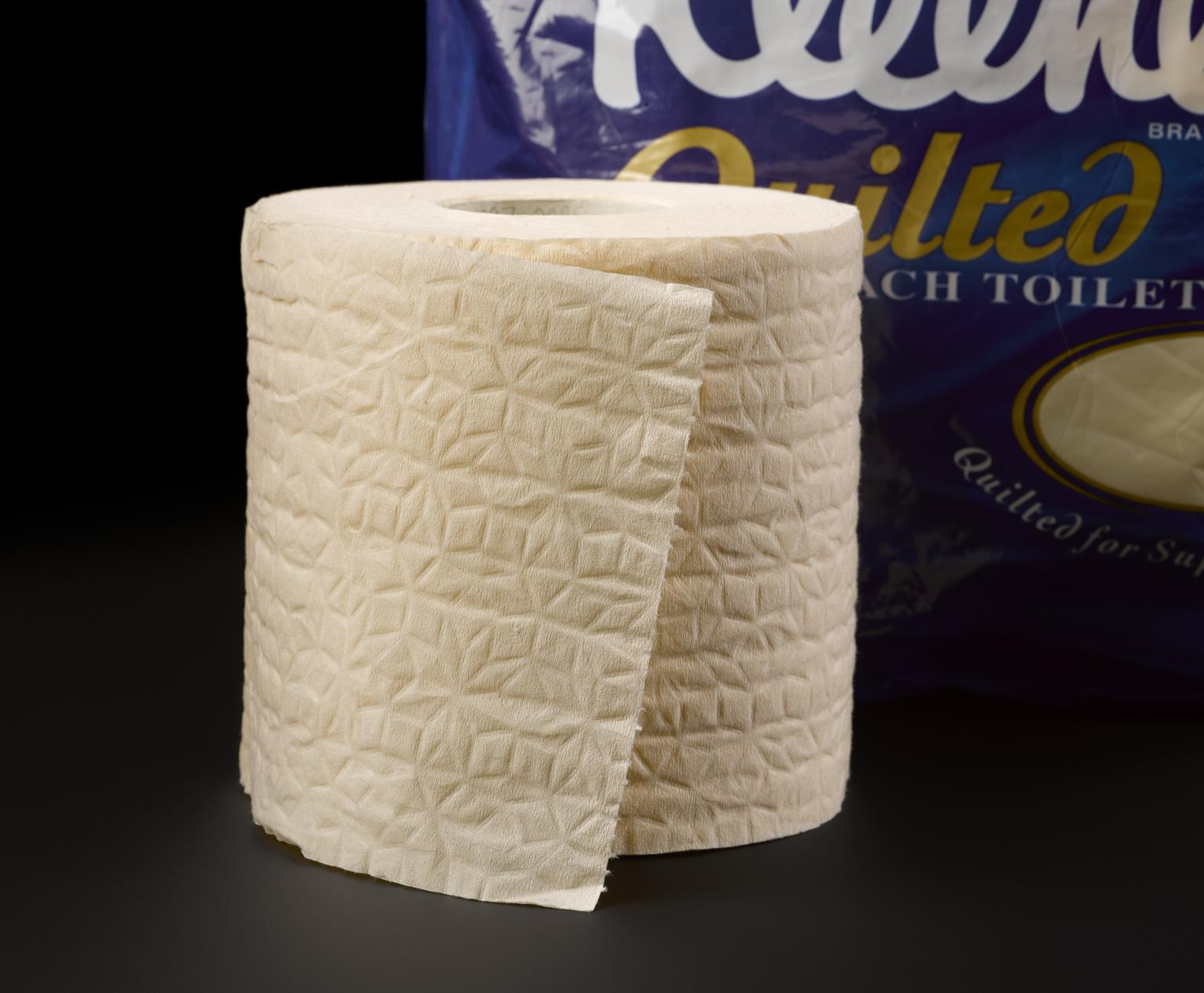}
\end{center}
\caption{Toilet paper seemingly embossed with Penrose tiles}
\label{Penrose_paper}
\end{figure}

The success of the Penrose tiles in popular culture has been immense, to the point that his tiles were seemingly imprinted on Kleenex toilet paper, see Figure \ref{Penrose_paper}. Luckily, the tiles had been patented, and the following was issued by David Bradley, the Director of Pentaplex:
\begin{quote}
So often we read of very large companies riding rough-shod over small businesses or individuals, but when it comes to the population of Great Britain being invited by a multi-national to wipe their bottoms on what appears to be the work of a Knight of the Realm without his permission, then a last stand must be made.
\end{quote}
The case was settled out of court and Kleenex stopped making the toilet paper. Squares are highly coveted by scientists studying aperiodic order, the second author has one framed on his office wall.

\section{The Taylor--Socolar monotile}

Given that two-tile aperiodic tile sets were found in the 1970s, it naturally begs the question of whether an aperiodic monotile is possible. That is, can one find a single tile and a set of local rules (decorations) that tile the plane but only nonperiodically. Here we have to be careful about what we mean by local rules. This is made very precise by Baake and Grimm in \cite[Section 5.7]{BG}. For us, suffice it to say that a local rule is a decoration that can only ``see" a finite and predefined radius in $\mathbb{R}^2$.

The Taylor--Socolar tile is an excellent example of a decoration that defines a local rule that is not edge-to-edge. The rules were discovered in 2010 by Joan Taylor, an amateur mathematician from Australia \cite{Tay}. She made contact with tiling expert Joshua Socolar in order to verify the discovery and fully work out the details. In a crowning achievement, they introduced the first true aperiodic monotile with local rules that are defined only between neighbours and next nearest neighbours \cite{ST, ST2}.

\begin{figure}[ht]
\begin{center}
\includegraphics[width=5.5cm]{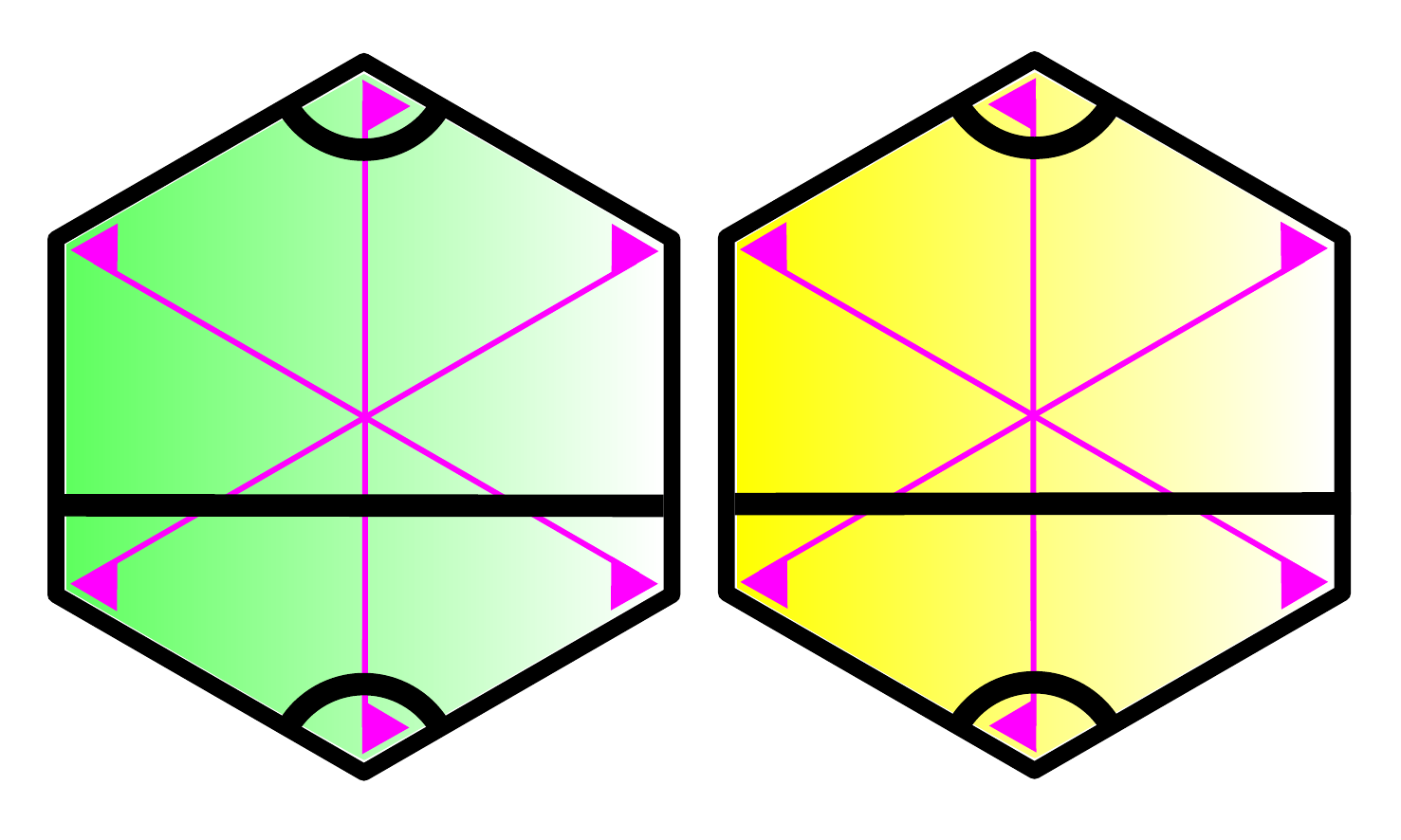} $\qquad$
\includegraphics[width=7cm]{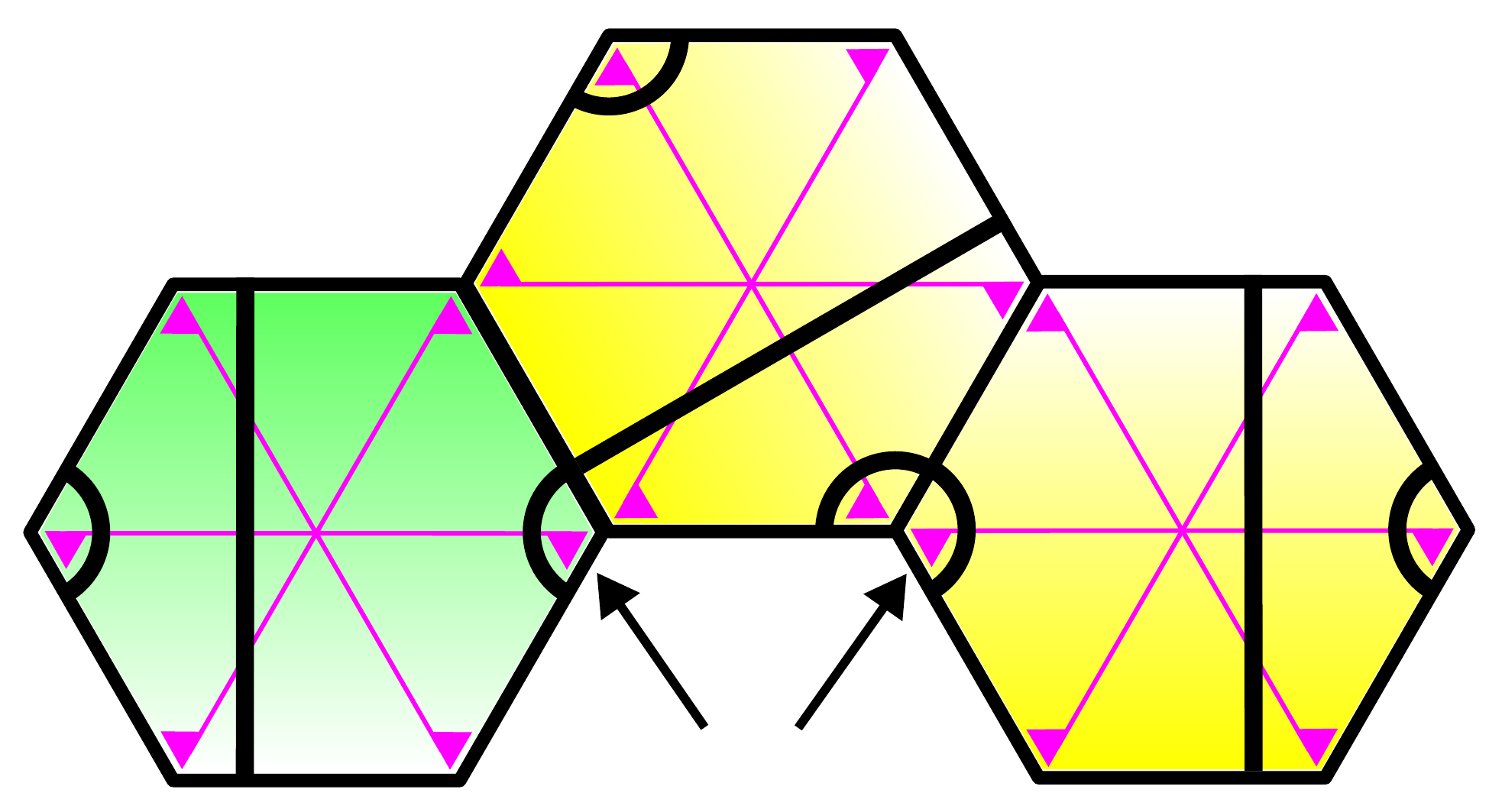}
\end{center}
\caption{The Taylor--Socolar monotile. On the left is the tile and its reflection. The image on the right demonstrates the (\textbf{R2})-rule, flags that meet across a single hexagonal edge must point in the same direction.}
\label{TS_tiles}
\end{figure}

\begin{figure}[ht]
\begin{center}
\includegraphics[width=15cm]{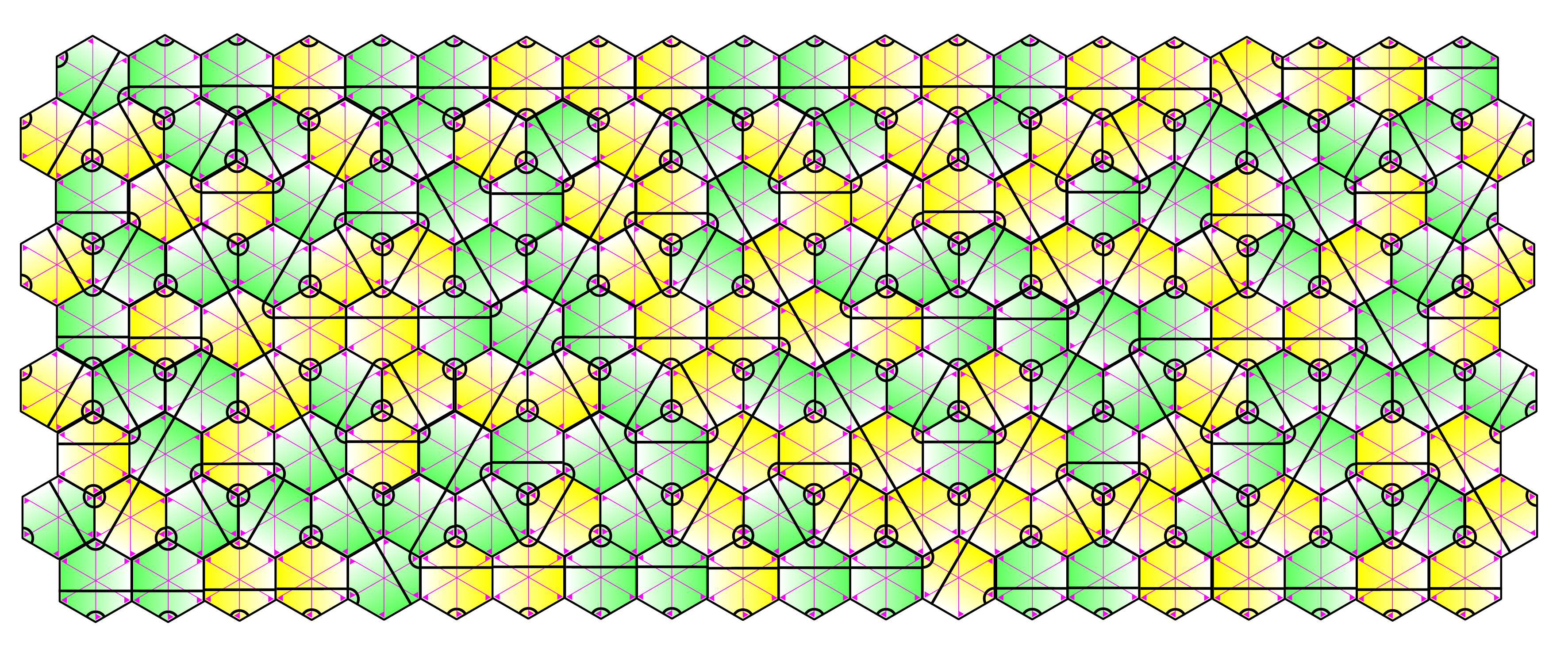}
\end{center}
\caption{A patch of a Taylor--Socolar tiling.}
\label{TS_patch}
\end{figure}

The Taylor--Socolar tile consists of a hexagonal tile and its mirror image such that:
\begin{enumerate}
\item[(\textbf{R1})] The black lines must continue across tile edges,
\item[(\textbf{R2})] The purple flags at vertices of tiles that meet across a single hexagonal edge must point in the same direction.
\end{enumerate}
The tiles are depicted in Figure \ref{TS_tiles} and the arrows on the right point to the flag rule (\textbf{R2}). Figure \ref{TS_patch} shows a patch and the reader should verify that (\textbf{R1}) and (\textbf{R2}) hold in any local patch in order to properly understand the local rules.

\begin{theorem}[{\cite[Theorem 1]{ST}}]\label{TS-theorem} The Taylor--Socolar monotile is aperiodic; that is, there are tilings formed by isometries of the Taylor--Socolar tile satisfying (\textbf{R1}) and (\textbf{R2}) in every local patch, and every such tiling is nonperiodic.
\end{theorem}

The proof of Theorem \ref{TS-theorem} is very elegant and is recommended to the interested reader.

\begin{figure}[ht]
\begin{center}
\includegraphics[width=7cm]{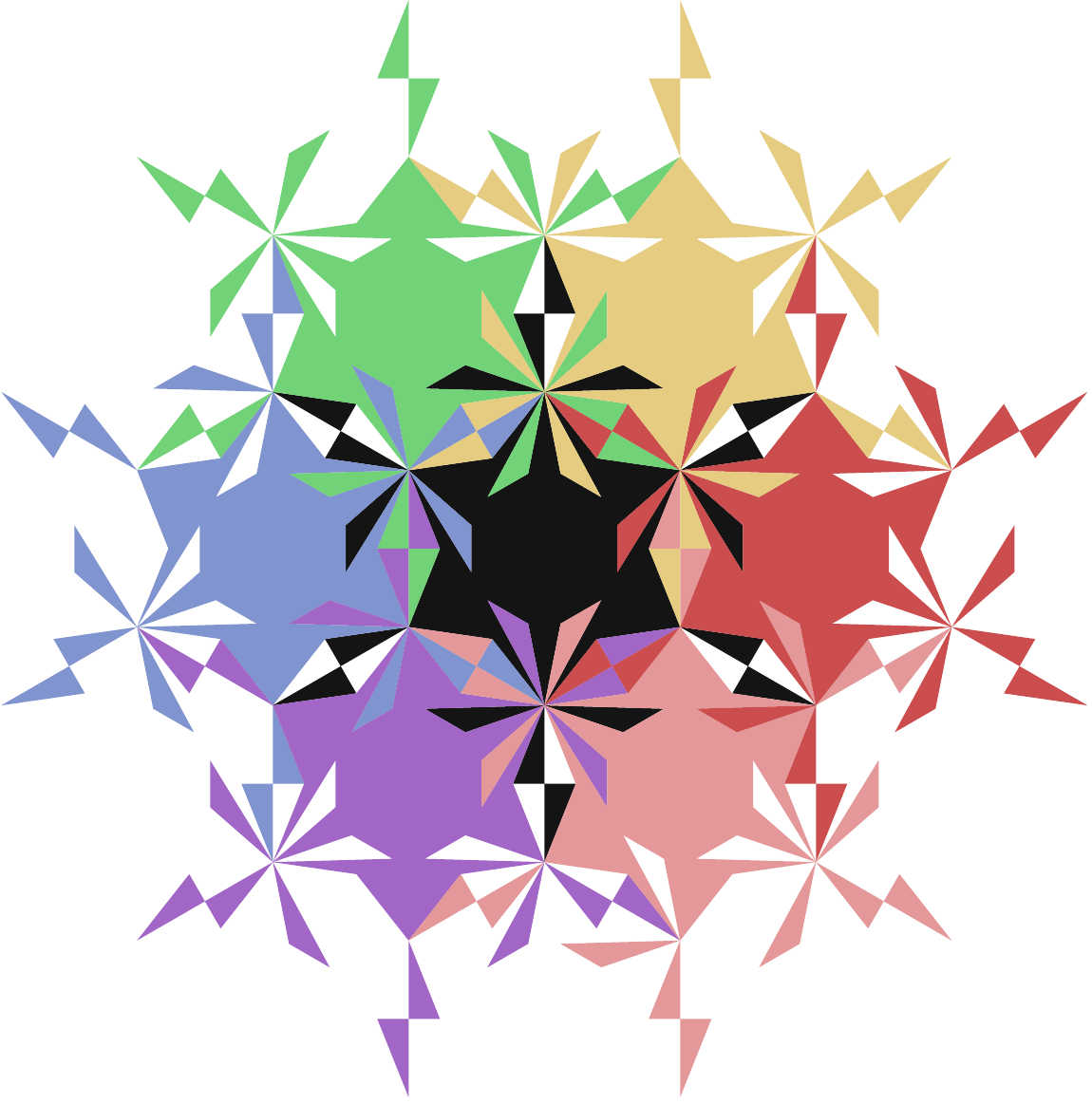}
\end{center}
\caption{A geometric version of the Taylor--Socolar tile. Distinct tiles have different colours so that we can see where the connected tiles interpenetrate. Image taken from \cite[Figure 6]{ST2}.}
\label{TS_geom}
\end{figure}

Interestingly, the Taylor--Socolar tile can be realised by shape alone. However, the tile is not simply connected; that is, the tile is connected but is not a topological disc. See Figure \ref{TS_geom}, where distinct tiles are coloured differently so that one can see how tiles interconnect. Focussing on the single black tile in the centre of the image shows the interesting geometry of the tile.

For the experts, adding an extra vertex consistency rule, (\textbf{R3}), to the Taylor--Socolar tile forces a single LI-class, see \cite[p.2215]{ST}. This additional rule was included in the original discovery by Taylor \cite{Tay}. In practice, this means that the tiling space of the Taylor--Socolar tiling satisfying (\textbf{R1})-(\textbf{R3}) is minimal in the sense of dynamical systems. Moreover, the Taylor--Socolar tilings are model sets \cite{BGG,LM}, and thus are strongly related to mathematical quasicrystals. These properties imply that the Taylor--Socolar tile is immensely important from the perspective of aperiodic order.

\begin{figure}[ht]
\begin{center}
\includegraphics[width=6cm]{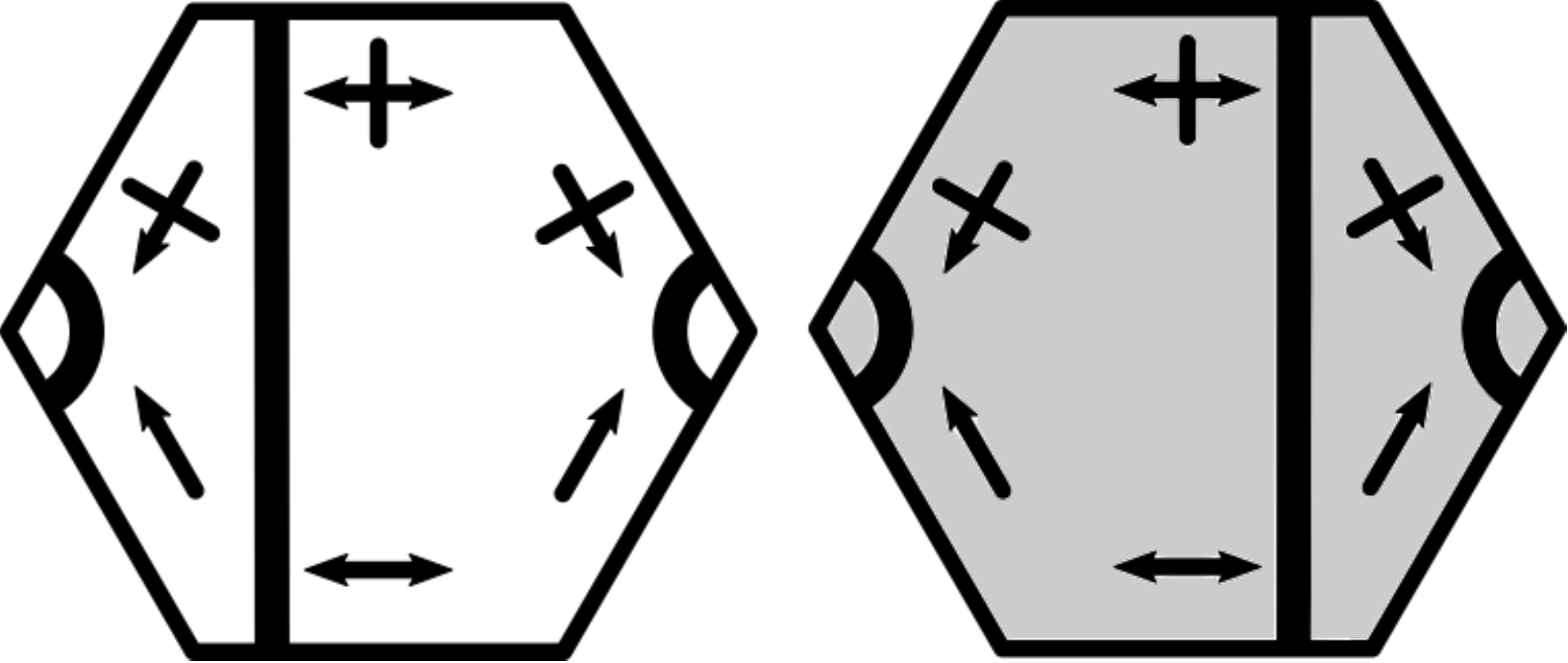}
\end{center}
\caption{An orientational monotile, where the local rule depends on the orientation of arrow in the plusses and minuses. The top edges are labelled with positive $\mathcal{R}2$-charges, which are oriented in each tile from left to right: anticlockwise, both clockwise and anticlockwise, and clockwise, respectively. The bottom edges are labelled with negative $\mathcal{R}2$-charges, from left to right: clockwise, both clockwise and anticlockwise, and anticlockwise}
\label{Orientational_tile}
\end{figure}

Jamie Walton and the second author found a modification of the Taylor--Socolar tile that has edge-to-edge matching rules \cite{WW}. However, the ($\mathcal{R}2$) rule below is somewhat unusual, it depends on orientation in the following way. Two tiles $t_1$ and $t_2$ are permitted to meet along a shared edge $e$ only if:
\begin{enumerate}
\item[($\mathcal{R}$1)] The black lines continue across $e$,
\item[($\mathcal{R}2$)] Whenever the two charges at $e$ in $t_1$ and $t_2$ both have a clockwise orientation then they must be opposite in charge.
\end{enumerate}

\begin{figure}[ht]
\begin{center}
\includegraphics[width=15cm]{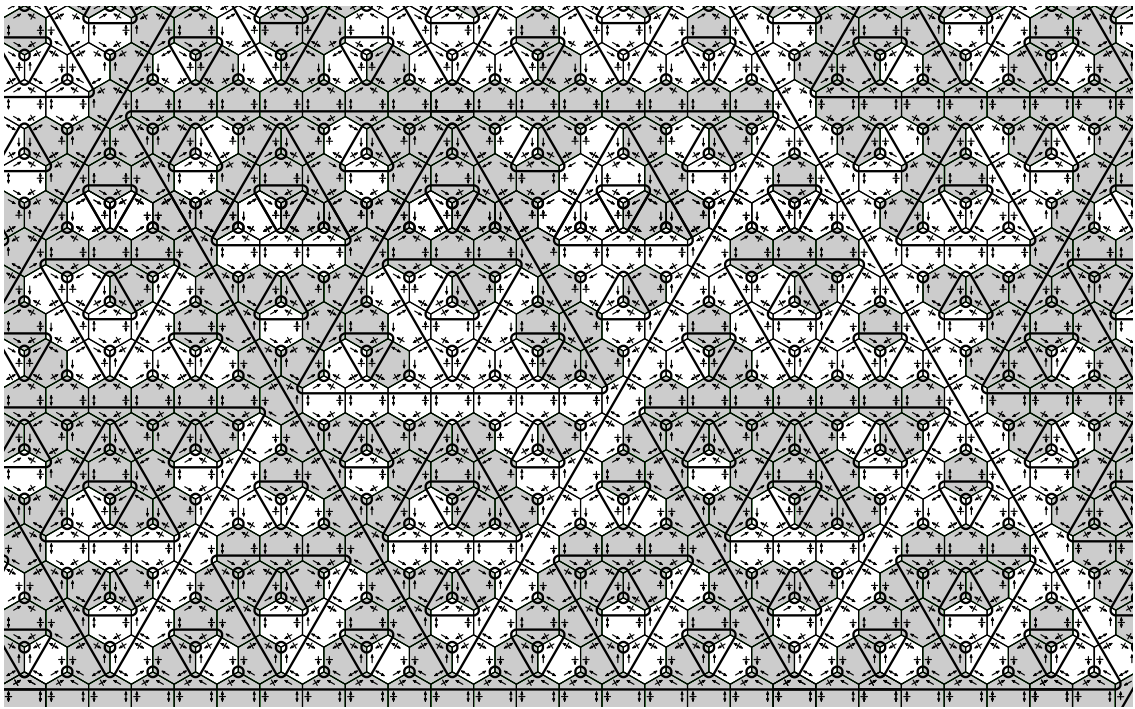}
\end{center}
\caption{A patch of the orientational monotile. Two neighbouring charges with clockwise orientation must be opposite.}
\label{Orientational_patch}
\end{figure}

\begin{theorem}[{\cite[Theorem 1.1]{WW}}]\label{Orientational-theorem} The tile in Figure \ref{Orientational_tile} is aperiodic; that is, there are tilings formed by isometries of the tile satisfying ($\mathcal{R}$1) and ($\mathcal{R}$2) in every local patch, and every such tiling is nonperiodic.
\end{theorem}

In this case the proof of aperiodicity is surprisingly simple and does not require a lot of case-checking. Indeed, we show that ($\mathcal{R}$1)-lines always lead to longer ($\mathcal{R}$1)-lines. Thus, there are arbitrarily long ($\mathcal{R}$1)-lines in any tiling. However, this precludes translational periodicity.  For if there was a non-trivial translation $x \in \mathbb{R}^2$ such that $T+x=T$, then there must be an ($\mathcal{R}$1)-line that is longer than $|x|$ and the structure of the ($\mathcal{R}$1)-triangles forbids this translation.

\section{The Hat and Spectre tiles}

The Hat monotile \cite{Hat} was discovered by David Smith, Craig S. Kaplan, Joseph Samuel Myers, and Chaim Goodman-Strauss with the article appearing on the Mathematics arXiv in March 2023. The paper generated an immediate buzz, resulting in newspaper articles in both the New York Times and the Guardian. The Hat was originally discovered by David Smith in November 2022 and the authors worked furiously to understand whether the Hat is an aperiodic monotile.

\begin{figure}[ht]
\begin{center}
\includegraphics[width=10cm]{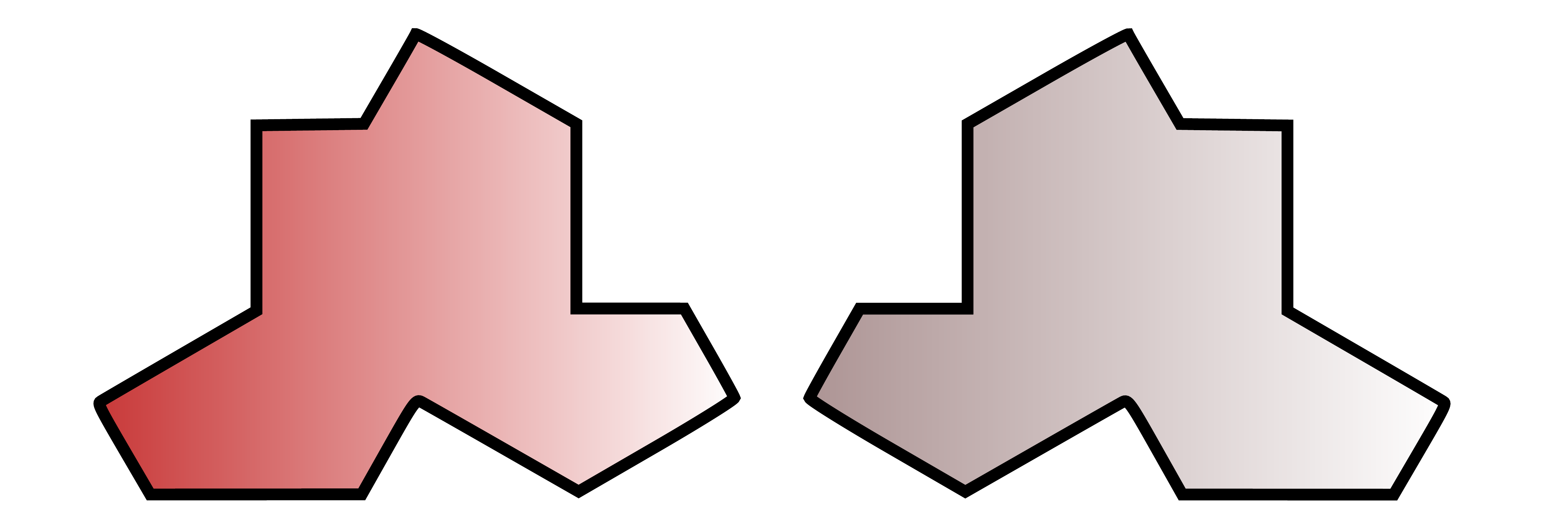} $\qquad$
\includegraphics[width=4.2cm]{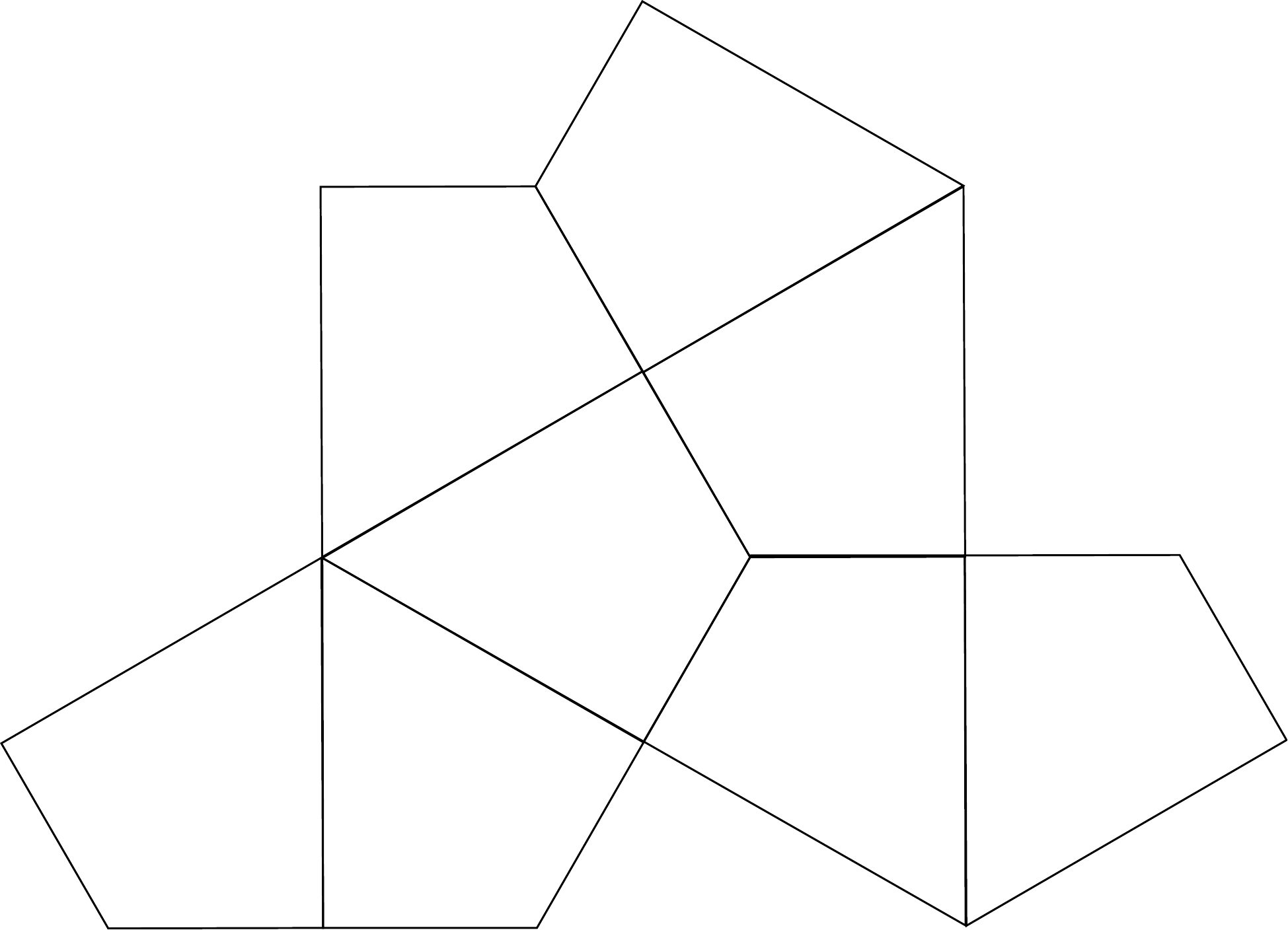}
\end{center}
\caption{The Hat tile and its mirror image are an aperiodic tile set without any need for further decorations. On the right we see how the Hat is formed from kites that combine to form regular hexagons.}
\label{HatTiles}
\end{figure}

The Hat tile is unbelievably simple and elegant. It can be found by forming a hexagonal grid, dividing each hexagon into kites and then combining kites from three neighbouring hexagons into a tile. The kites are formed by cutting the hexagons with straight lines through the midpoints of opposite edges. For this reason the authors often refer to it as a polykite. See Figure \ref{HatTiles} for the Hat tile, its mirror image, and a rendering into kites that combine to form regular hexagons.

Interestingly, the Hat was discovered by Smith when tinkering with polyforms to see what sorts of visually interesting tiling he could create. He got stuck building large patches of Hat tiles and sought help from Kaplan's Heesch Number software. The \emph{Heesch Number} of a tile is the largest number of concentric rings that  form a patch around a single tile by isometric copies, where a ring consists of all tiles that touch the previous ring. The current record holder for a tile that does not tile the plane was found by Ba\v{s}i\'c \cite{Bas} and has Heesch Number 6. Kaplan was able to show that the Hat has at least Heesch number 10 and then improved that to 16. So, it seemed pretty likely that the Hat would tile the plane, and in a way that didn't seem to have any translational periodicity. This type of attack on the monotile problem was completely different from previous attempts, where authors typically start with a tiling of the plane and then add rules to try and enforce nonperiodicity. In some ways this approach goes back to a prototile set decidability problem, in the original sense of Wang.

\begin{figure}[ht]
\begin{center}
\includegraphics[width=15cm]{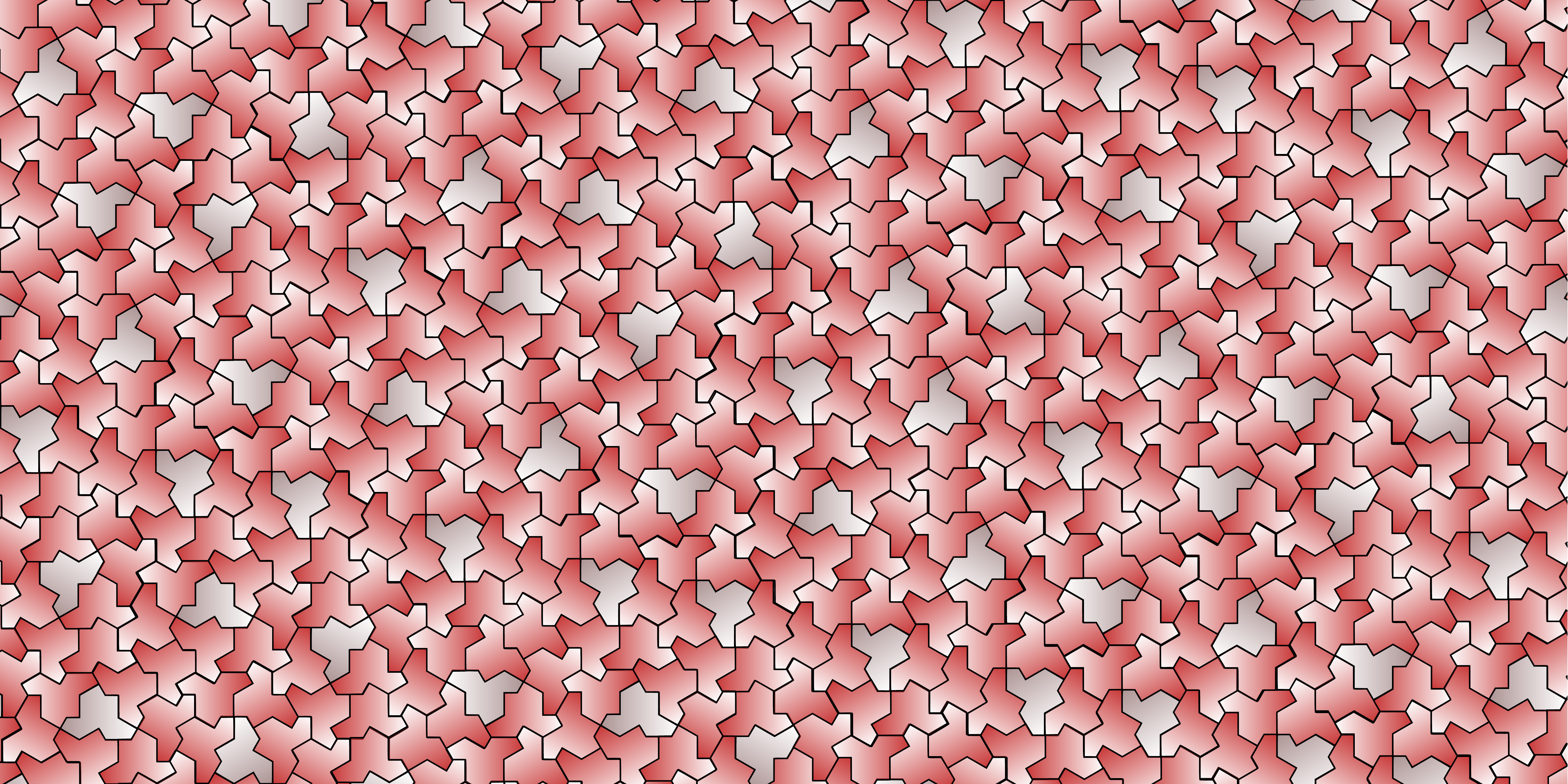} 
\end{center}
\caption{A patch of Hat tiles.}
\label{HatPatch}
\end{figure}

A patch of an infinite tiling formed from Hats is depicted in Figure \ref{HatPatch}. Notice that reflected Hats appear with a low frequency compared to Hat tiles.

\begin{theorem}[{\cite[Theorem 1.1]{Hat}}]\label{Hat-theorem} The Hat monotile is aperiodic; that is, there are tilings formed by isometries of the Hat tile, and every such tiling is nonperiodic.
\end{theorem}

There are two proofs that the Hat tiles the plane and two proofs of aperiodicity. The authors show that Hat tilings arise from a (nonstone) substitution rule that is recognisable, in the sense that the hierarchy can be deduced in an infinite tiling. On one hand, this allows one to construct patches of arbitrarily large size that nest into each other in order to construct a tiling of the plane. One the other hand, recognisability implies that one can identify structure in the tiling of arbitrarily large size, and hence the tiling cannot be periodic. The second proof of existence uses a rather simple fusion system (see \cite{FS} for the definition of fusion) to define arbitrarily large patches of tiles that expand out from a fixed tile, see \cite[Figure 2.11]{Hat}. The second proof of aperiodicity is more involved but has already led to the discovery of the Spectre tile.

The idea behind the second proof of aperiodicity is to contract the edges of the hat tile to find other combinatorially equivalent tilings in the sense that the patterns formed by tiles are the same. Indeed, the two edge lengths in a Hat tile are at a ratio of $\sqrt{3}$ to one another and come in complementary (opposite) pairs. Thus, we can label the Hat tile as Tile$(\sqrt{3},1)$. The idea is now to consider Tile$(a,b)$ for $0 \leq a \leq \sqrt{3}$ and $0 \leq b \leq 1$. There is an excellent animation of this on YouTube titled "Aperiodic monotile animation".  At the two extremes are Tile$(0,1)$, called the Comet and Tile$(\sqrt{3},0)$ called the Chevron, see \cite[Figure 3.1]{Hat}. Since these two tiles are combinatorially equivalent to the Hat tiling, the authors are free to use them to deduce properties of the Hat tiling. Thus, the authors suppose that the Comet tiling is periodic and aim to derive a contradiction, meaning that this hypothesis could not be correct and the tiling is nonperiodic. Since the Comet is assumed to be periodic and forms the same combinatorial tiling as the Chevron, we can deduce that the Chevron is also periodic. Moreover, since the tilings are combinatorially equivalent there must be an affine map between the periodic lattice of the Comet and that of the Chevron. Using an argument, somewhat similar to an argument that $\sqrt{2}$ is irrational, they prove that such an affine map cannot exist. This implies that the Hat tiling must also be nonperiodic. See \cite[Section 3]{Hat} for further details.

One of the most interesting developments was again discovered by Smith and his coauthors, the Spectre tile \cite{Spectre}. Amazingly, this is Tile$(1,1)$ from the previous paragraph, which can be used to tile the plane periodically if one allows reflection of the tile. However, the authors realised that it is still possible to tile the plane if reflections of the tile are not allowed to appear in a tiling, and more amazingly that all such tilings are nonperiodic!

\begin{figure}[ht]
\begin{center}
\includegraphics[width=10cm]{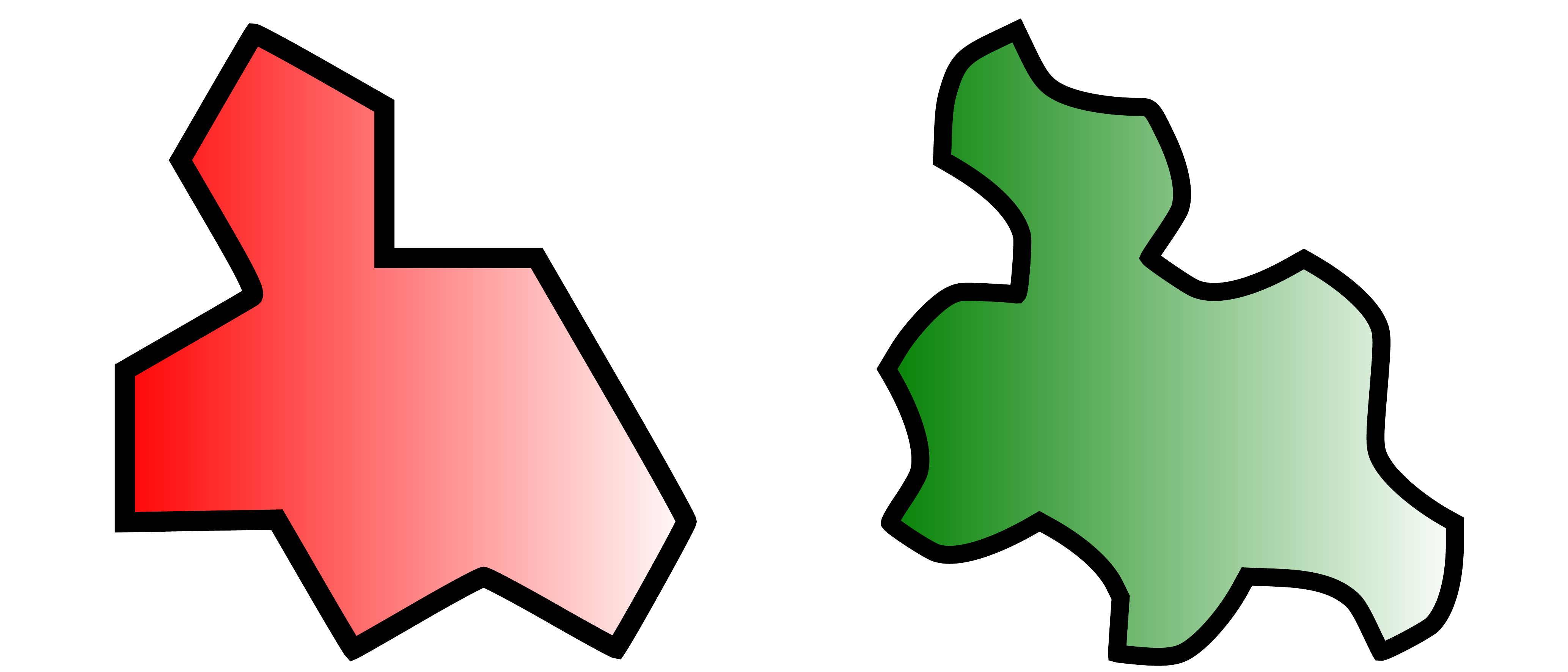} $\qquad$
\includegraphics[width=4cm]{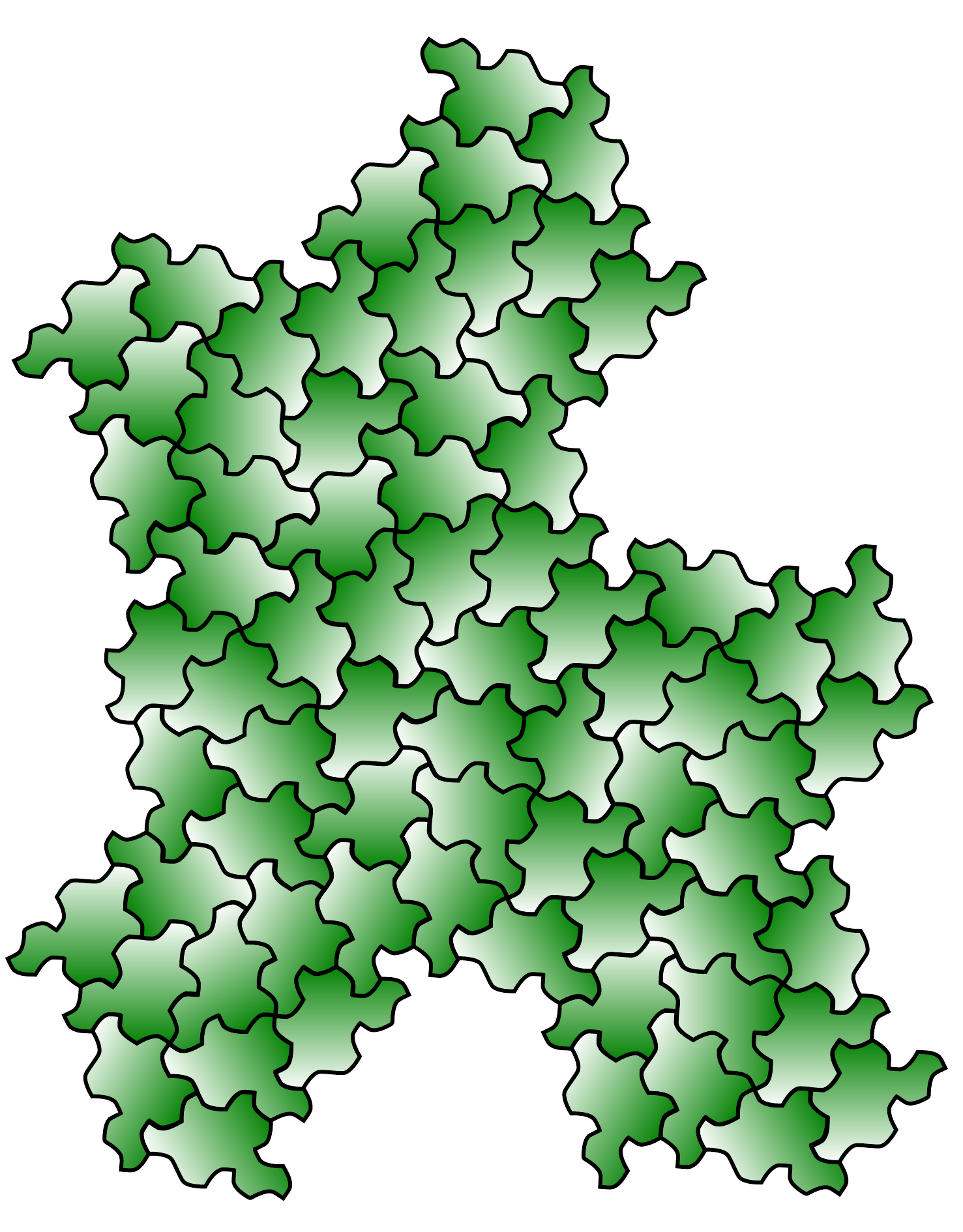}
\end{center}
\caption{The image on the left is Tile$(1,1)$ and one must forbid reflections to obtain an aperiodic monotile. The tile in the centre is the Spectre tile that tiles the plane without allowing reflection, and does not tile the plane if there is at least one Spectre tile and at least one mirror image of the Spectre tile. The image on the right comes from a fusion system to construct arbitrarily large patches of Spectre tiles.}
\label{SpectreTiles}
\end{figure}

\begin{figure}[ht]
\begin{center}
\includegraphics[width=15cm]{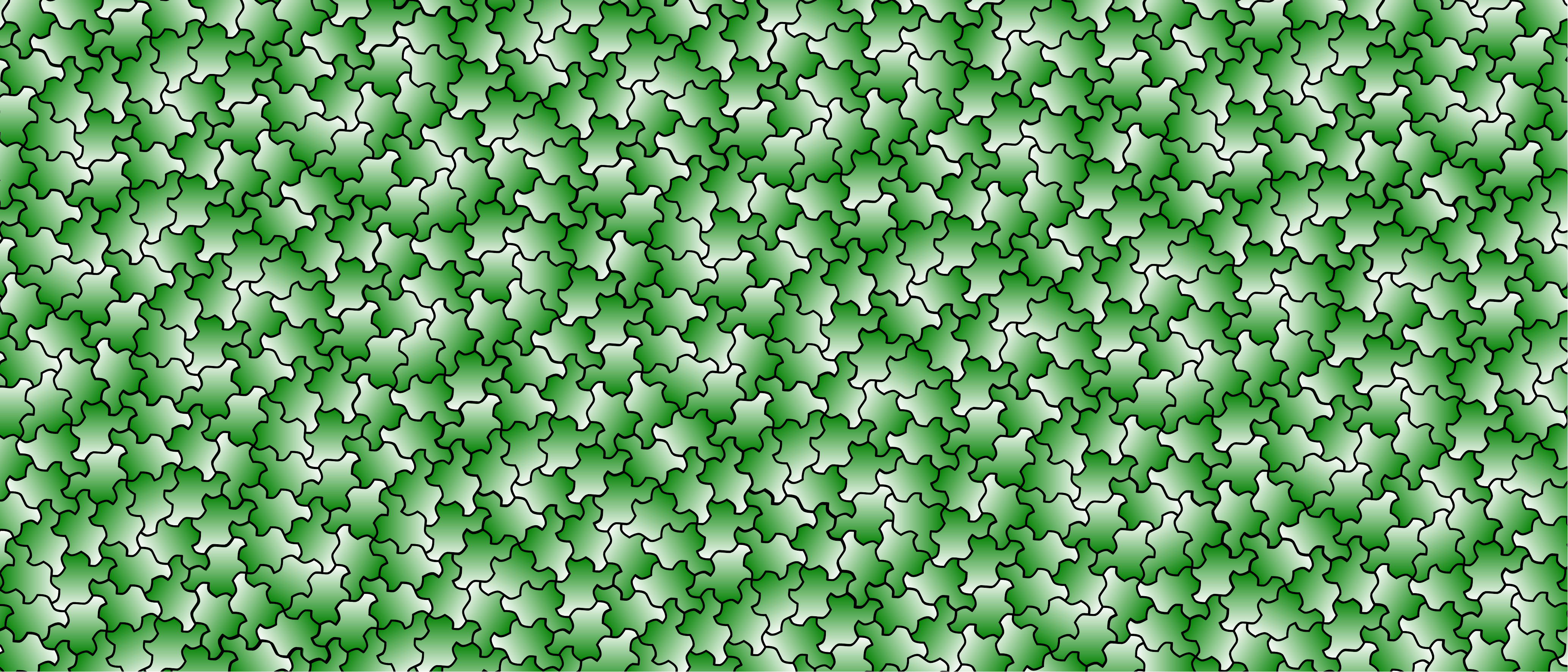}
\end{center}
\caption{A patch of Spectre tiles.}
\label{Spectre}
\end{figure}

Tile$(1,1)$ appears on the left hand side of Figure \ref{SpectreTiles} and a version of the Spectre tile appears in the centre. The Spectre is merely an edge modification of Tile$(1,1)$ to curves, which eliminate the possibility of using both the tile and its reflection to tile the plane. We understand that Dave Smith proposed to call the tile the \emph{Spectre} due to the image in the centre of Figure \ref{SpectreTiles}. The image on the right was constructed through a fusion system to build arbitrarily large patches \cite{Spectre}, and we think that also looks spectre like.

\begin{theorem}[{\cite[Theorem 2.2]{Spectre}}]\label{Spectre-theorem} Tile$(1,1)$ and the Spectre monotile are aperiodic; that is, there are tilings formed by direct (orientation-preserving) isometries of each tile, and every such tiling is nonperiodic.
\end{theorem}

Let us make a couple of remarks about Theorem \ref{Spectre-theorem}. First, this is an absolutely incredible result that was completely unexpected, even given the recent Hat tile result. We note that Tile$(1,1)$ is referred to by the authors of \cite{Spectre} as a \emph{weakly chiral monotile} since it satisfies Theorem \ref{Spectre-theorem}, but allowing a reflection results in a tile that can be used to construct a periodic tiling. A Spectre tile, see the centre of Figure \ref{SpectreTiles}, is referred to as a \emph{strictly chiral aperiodic monotile} since it satisfies Theorem \ref{Spectre-theorem}, but any prototile set containing both the Spectre and its mirror image nonredundantly does not tile the plane.

Theorem \ref{Spectre-theorem} was proved by showing that Tile$(1,1)$ and the Spectre arise from a (nonstone) inflation rule that is recognisable, similarly to the first proof of Theorem \ref{Hat-theorem} for the Hat tile.

For the experts, Baake, G\"ahler, and Sadun have extended the 1-parameter family of tiles given by Tile$(a,b)$ to complex variables \cite{BGS}. They showed that all the continuous hulls are topologically conjugate dynamical systems under these parameters, up to linear rescaling of the ambient space, and found a self-similar representative they call the CAP tiling. The name follows from their result that the tiling is a \emph{c}ut \emph{a}nd \emph{p}roject tiling, and hence forms a model set. They also compute the cohomology of the family and show that it has pure-point dynamical spectrum. In current work, the same authors have also shown that the Spectre tile has similar properties \cite{BGS2}.

Since the original Hat preprint appeared \cite{Hat}, Akiyama and Araki have provided yet another proof of existence and aperiodicty of a member of the Hat family \cite{AA}. They use the Golden Hex substitution to prove existence and Golden Ammann bars to prove aperiodicty.

Mathematicians are now discovering that the unique structure of the Spectre tiling lends it many other fascinating properties \cite{Felix}. For example, the \emph{dimer model} asks how many ways there are to colour the edges of the tiles so that each vertex meets exactly one coloured edge (dimer). Remarkably, this model can be exactly solved on the Spectre tiling. The number of dimer arrangements is $2^{N_{\textrm{Mystic}}+1}$, where $N_{\textrm{Mystic}}$ is the number of \emph{Mystic} tiles, see \cite[p.6]{Spectre} or \cite[p.2]{Felix} for the definition of a Mystic. More remarkable still, the dimer model can also be exactly solved when quantum superpositions of dimer placements are allowed! Thus, Singh and Flicker have exactly solved the quantum dimer model for the first time in any setting.

\end{document}